\documentclass[11pt,reqno]{amsart}
\usepackage{amssymb}

\newcommand{\map}[1]{\;\xrightarrow{#1}\;}
\newcommand{\iso}{\cong}
\newcommand{\mil}{\varprojlim}
\newcommand{\dlim}{\varinjlim}

\newcommand{\fa}{\mathfrak a}
\newcommand{\ep}{\varepsilon}
\newcommand{\cK}{\mathcal K}
\newcommand{\cL}{\mathcal L}
\newcommand{\Q}{\mathbf Q}
\newcommand{\bR}{\mathbf R}
\newcommand{\Z}{\mathbf Z}
\newcommand{\cO}{\mathcal O}
\newcommand{\co}{\mathcal O}
\newcommand{\p}{\mathfrak p}
\newcommand{\fp}{\mathfrak p}
\newcommand{\bE}{\mathbf E}
\newcommand{\ff}{\mathfrak f}
\newcommand{\Ell}{\mathrm {ell}}
\newcommand{\tors}{\mathrm {tor}}

\newcommand{\Sel}{\mathrm {Sel}}
\newcommand{\str}{\mathrm {str}}
\newcommand{\rel}{\mathrm {rel}}
\newcommand{\CS}{\mathcal S}
\newcommand{\cS}{\mathcal S}
\newcommand{\HH}{\mathcal H}
\newcommand{\cH}{\mathcal H}
\newcommand{\CC}{\mathcal C}
\newcommand{\cC}{\mathcal C}

\newcommand{\bQ}{\mathbf Q}
\newcommand{\bZ}{\mathbf Z}

\DeclareMathOperator{\Aut}{Aut}
\DeclareMathOperator{\cha}{Char}

\DeclareMathOperator{\Hom}{\mathrm{Hom}}
\DeclareMathOperator{\Gal}{\mathrm{Gal}}
\DeclareMathOperator{\loc}{\mathrm{loc}}

\DeclareMathOperator{\rk}{\mathrm {rk}}

\input xy
\xyoption{all}

\begin{document}
\title{Anticyclotomic Iwasawa theory of CM elliptic curves II}
\author{Adebisi Agboola \and Benjamin Howard}
\date{Final version, March 23, 2005}

\address{Department of Mathematics\\ University of California\\ Santa
  Barbara, CA\\ 93106}
\email{agboola@math.ucsb.edu}
\address{Department of Mathematics\\ Harvard University\\ Cambridge, MA\\
02138}
\curraddr{Department of Mathematics\\ University of Chicago\\ 
\\ 5734 S. University Ave.\\ Chicago, IL\\
60637}
\email{howard@math.uchicago.edu}
\thanks{The first author was partially supported by NSF
grant DMS-0070449}
\thanks{The second author was supported by an NSF
Mathematical Sciences Postdoctoral Research Fellowship.}

\begin{abstract}
We study the Iwasawa theory of a CM elliptic curve $E$ in the
anticyclotomic $\mathbf{Z}_p$-extension $D_\infty$ of the CM field $K$,
where $p$ is a prime of good, supersingular reduction for $E$. Our
main result yields an asymptotic formula for the corank of the
$p$-primary Selmer group of $E$ along the extension $D_\infty/K$.
\end{abstract}

\maketitle

\theoremstyle{plain}
\newtheorem{BigThm}{Theorem}
\newtheorem{Thm}{Theorem}[section]
\newtheorem{Prop}[Thm]{Proposition}
\newtheorem{Lem}[Thm]{Lemma}
\newtheorem{Cor}[Thm]{Corollary}
\newtheorem{Conj}[Thm]{Conjecture}
\newtheorem{Assum}[Thm]{Assumption}

\theoremstyle{definition}
\newtheorem{Def}[Thm]{Definition}
\newtheorem{Rem}[Thm]{Remark}
\newtheorem{Ques}[Thm]{Question}

\renewcommand{\theBigThm}{\Alph{BigThm}}
\numberwithin{equation}{section}

\setcounter{section}{0}


\section{Introduction} \label{S:intro}


Let $E$ be an elliptic curve over $\Q$. Whilst much is known about the
Iwasawa theory of $E$ for primes of ordinary reduction, the same is
unfortunately not true of Iwasawa theory at supersingular primes, for
in this case the Iwasawa modules that one naturally considers are not
torsion, and the obvious candidates for $p$-adic $L$-functions do not
lie in the Iwasawa algebra. Nevertheless, there has recently been a
great deal of progress in the study of the Iwasawa theory of elliptic
curves at supersingular primes. In particular, S. Kobayashi has
recently formulated a cyclotomic main conjecture for $E$ within this
framework (see \cite{Ko}). His conjecture relates certain restricted
`plus/minus' Selmer groups of $E$ to certain modified $p$-adic
$L$-functions defined by R. Pollack (see \cite{Po}), and it is
equivalent to a cyclotomic main conjecture that was proposed earlier
by K. Kato and B. Perrin-Riou (see \cite{Ka}, \cite{PR5}
\cite{PR4}). Kobayashi's conjecture has recently been proved by Rubin
and Pollack (see \cite{RP}) when $E$ has complex multiplication, and
Kobayashi himself, using methods of Kato, proves one divisibility of
the main conjecture in the non-CM case. In both cases, the plus/minus
Selmer groups are cotorsion modules over the cylotomic Iwasawa
algebra, and so the corank of the $p$-Selmer group remains bounded as
one ascends the cyclotomic $\Z_p$-extension.

Suppose now that $E$ has complex multiplication by the maximal order
$\co$ of an imaginary quadratic field $K$.  Let $\psi$ denote the
$K$-valued grossencharacter associated to $E$, and write $\ff$ for the
conductor of $\psi$.  Fix once and for all a rational prime $p>3$ at
which $E$ has good reduction, and which is inert in $K$.  Then $E$ has
supersingular reduction at $p$.  We write $\p$ for the unique prime of
$K$ above $p$.  Let $D_\infty$ be the anticyclotomic $\Z_p$ extension
of $K$, and let $D_n\subset D_\infty$ be the subfield such that
$[D_n:K]=p^n$.  The prime $\p$ is totally ramified in $D_\infty$, and
we let $\p$ also denote the unique place of $D_\infty$ above $\p$.

In this paper, we study the Iwasawa theory of $E$ over $D_\infty$. We
define anticyclotomic versions of Kobayashi's restricted plus/minus
Selmer groups, and we analyse their structure using the Euler system
of twisted elliptic units (cf. \cite{AH}). In contrast to what happens
in the cyclotomic case, it turns out that one of the restricted Selmer
groups is a cotorsion Iwasawa module, while the other is not, and
which module is cotorsion is determined by the sign in the functional
equation of $L(E/\Q,s)$. Our main result, predicted by R. Greenberg
\cite[p. 247]{Gr}, is as follows (somewhat more information is
contained in Theorem \ref{final}).

\begin{BigThm}
Let $\phi$ be Euler's function and let $\ep=\pm 1$ be the sign in the
functional equation of $L(E/\Q,s)$. Write $\Sel_{\p^\infty}(E/D_n)$
for the $\p$-primary Selmer group of $E/D_n$, and $\co_\p$ for the
local completion of $\co$ at $\p$. Then there is an integer $e$,
independent of $n$, such that
$$
\mathrm{corank}_{\co_\p}\ \Sel_{\p^\infty}(E/D_n)=
e+\sum_{1\le k\le n,\ (-1)^k=\ep}\phi(p^k)
$$
for all $n\gg 0$.
\end{BigThm}

The results in this paper may be viewed as a first step towards a
supersingular main conjecture of the same type as that considered in
\cite{AH}.  However, in the present setting, we do not know how to
define suitable anticyclotomic analogues of Pollack's plus/minus
$p$-adic $L$-functions.  The essential missing ingredient is a
construction of local elements along the lines of \cite[\S 8.4]{Ko}.


\section{Selmer groups} \label{S:selmer}


We write
$$
T=T_p(E),\quad W=E[p^\infty]
$$
for the $p$-adic Tate module and the goup of $p$-power torsion points
in $E(\overline{K})$ respectively.  Let $F/K$ be any finite
extension. For any place $v$ of $F$, we define $H^1_f(F_v,W)$ to be
the image of $E(F_v) \otimes (\Q_p/\Z_p)$ under the Kummer map
$$
E(F_v) \otimes (\Q_p/\Z_p) \to H^1(F_v,W),
$$
and we write $H^1_f(F_v,T)$ for the orthogonal complement of
$H^1_f(F_v,W)$ with respect to the local Tate pairing. 
Note that $H^1_f(F_v,W)=0$ if $v\nmid p$.
If $c \in H^1(F,W)$, then we write $\loc_v(c)$ for the image of $c$ in
$H^1(F_v,W)$. 

We define

$\bullet$ the {\it relaxed Selmer group} $\Sel_{\rel}(F,W)$ by
$$
\Sel_{\rel}(F,W) = \left\{ c \in H^1(F,W) \mid \loc_v(c) \in
H^1_f(F_v,W)\, \text{for all $v$ not dividing $p$}
\right\};
$$

$\bullet$ the {\it true Selmer group} $\Sel(F,W)$ by
$$
\Sel(F,W) = \left\{ c \in H^1(F,W) \mid \loc_v(c) \in
H^1_f(F_v,W)\, \text{for all $v$} \right\};
$$

$\bullet$ the {\it strict Selmer group} $\Sel_{\str}(F,W)$ by
$$
\Sel_{\str}(F,W) = \left\{ c \in \Sel(F,W) \mid \loc_v(c) = 0\,
\text{for all $v$ dividing $p$} \right\}.
$$

We also define $\Sel_{\rel}(F,T)$, $\Sel(F,T)$ and $\Sel_{\str}(F,T)$
in a similar way. It follows from the definitions that there are
inclusions
$$
\Sel_{\str}(F,W) \subset \Sel(F,W) \subset \Sel_{\rel}(F,W),
$$
and similarly with $W$ replaced by $T$.
If $F/K$ is an infinite extension, we define
$$
\Sel_*(F,W) = \varinjlim \Sel_*(F',W)
\hspace{1cm}
\cS_*(F,T) = \varprojlim \Sel_*(F',T),
$$
where the limits are taken with respect to restriction and
corestriction, respectively, over all subfields $F' \subset F$ finite
over $K$.

We now give the definition of a slightly modified (see Remark
\ref{R:trivial} below) form of Kobayashi's restricted plus/minus
Selmer groups. Let $\bE$ denote the formal group of $E$ over
$K_\fp$. Since $E$ has supersingular reduction at $p>3$,
it is a standard fact that $\bE$ is isomorphic to the unique
(up to isomorphism) Lubin-Tate formal goup over $K_\fp$ with 
parameter $-p$.  
For $n\ge 0$, let $\Xi^-_n$ be the set of characters of $\Gamma_n$ of exact
order $p^k$ with $k$ odd, together with the trivial character.  Let
$\Xi^+_n$ be the set of characters of $\Gamma_n$ of exact order $p^k$
with $k$ even, excluding the trivial character.  Define subspaces of
$E(D_{n,\p})\otimes_{\co_\p}K_\p$ by
$$
E_\pm(D_{n,\p})=\left\{ x\in \bE(D_{n,\p})\otimes K_\p\left|
\sum_{\sigma\in\Gamma_n}\chi(\sigma)x^\sigma=0,\ \forall \chi\in \Xi^\mp_n
\right.\right\}.
$$
Let $H^1_\pm(D_{n,\p},W)$ be the image of $E_\pm(D_{n,\p})$
under the Kummer map
$$
E(D_{n,\p})\otimes K_\p\map{} E(D_{n,\p})\otimes (K_\p/\co_\p)
\map{}H^1(D_{n,\p},W)
$$ 
and let $H^1_\pm(D_{n,\p},T)$ be the orthogonal complement of
$H^1_\pm(D_{n,\p},W)$ with respect to the local Tate pairing. We
define
$$
\Sel_\pm(D_n,W) = \left\{ c \in \Sel_{\rel}(D_n,W) \mid \loc_\fp(c)
\in H^1_\pm(D_{n,\fp},W) \right\};
$$
$$
\Sel_\pm(D_n,T) = \left\{ c \in \Sel_{\rel}(D_n,T) \mid \loc_\fp(c)
\in H^1_\pm(D_{n,\fp},T) \right\}.
$$
It follows from the definitions that we have inclusions
$$
\Sel_{\str}(D_n,W) \subset \Sel_{\pm}(D_n,W) \subset \Sel(D_n,W),
$$
$$
\Sel(D_n,T) \subset \Sel_{\pm}(D_n,T) \subset \Sel_\rel(D_n,T).
$$
In the limit, we define
$$
\Sel_\pm(D_\infty,W) = \varinjlim \Sel_\pm(D_n,W)
\hspace{1cm}
\cS_\pm(D_\infty,T) = \varprojlim \Sel_\pm(D_n,T),
$$
where the inverse limits are taken with respect to restriction and
corestriction, respectively.

In order to ease notation, we shall sometimes write
$$
\Sel^{\infty}_{*} = \Sel_*(D_\infty,W)\hspace{1cm} \cS_* =
\cS_*(D_{\infty},T).
$$


\section{Ranks} \label{S:ranks}


Let $\fa$ be an integral ideal of $\co$ coprime to $6p\ff$, and write
$\cK_{\fa}$ for the union of all ray class fields of $K$ of conductor
prime to $\fa$. Let $c_{\Ell,\fa}$ denote the Euler system of elliptic
units for $(\Z_p(1),\ff p,\cK_{\fa})$ in the sense of
\cite{Ru}. Twisting $c_{\Ell,\fa}$ by the character
$\omega_\mathrm{cyc}^{-1}\psi_{\p}$, where
$$
\psi_\p: \Gal(\overline{K}/K) \map{} \Aut_{\co_K}(T)\iso\co_\p^\times,
$$
and $\omega_{\mathrm{cyc}}$ is the cyclotomic character,
yields an Euler system $c_{\fa}$ for $(T,\ff p,\cK_{\fa})$ (see
\cite[Chapter 6]{Ru}). Then $c_{\fa}(F) \in \Sel_{\rel}(F,T)$ for
every $F \subset \cK_\fa$ finite over $K$. For any $L \subset
K_\infty$, let
$$
c_\fa(L) = \varprojlim c_\fa(F') \in \cS_\rel (L,T),
$$
where the inverse limit is taken over all subfields $F'$ of $L$ that
are finite over $K$. Let $\cC_\fa(F)$ be the
$\cO_{\fp}[[\Gal(F/K]]$-submodule of $\cS_\rel (F,T)$ generated by
$c_\fa (F)$, and write $\cC(F)$ for the submodule generated by
$\cC_\fa(F)$ as $\fa$ varies over all ideals that are coprime to
$6p\ff$. We set $\cC = \cC(D_\infty)$.

Define 
$$
\HH^1_\pm=\mil H^1_\pm(D_{n,\p},T)\hspace{1cm}
\HH^1=\mil H^1(D_{n,\p},T)
$$
and
$$
H^1_\pm(D_{\infty,\p},W)=\dlim H^1_\pm(D_{n,\p},W).
$$


Let $W(\psi)$ denote the root number of $\psi$. In particular
$W(\psi) = \pm 1$ and is equal to the sign in the functional equation 
of $L(E/\Q,s)$.

\begin{Prop} \label{P:nontriv}
The image of $\CC$ in $\HH^1$ is nontrivial, and lies in
$\HH^1_\ep$ if and only if $\ep$ is equal to the
sign of $W(\psi)$.
In particular $\CC$ is nontrivial, and is contained
in $\CS_\ep$ if and only if $\ep$ is the sign of $W(\psi)$.
\end{Prop}

\begin{proof} 
Let $\chi$ be any primitive character of $\Gal(D_n/K)$ for $n>0$, and
write $W(\chi \psi)$ for the root number of $\chi \psi$. The following
formula is proved by Greenberg in \cite[page 247]{Gr}:
\begin{equation} \label{E:rootno}
W(\chi\psi)=(-1)^{n+1}W(\psi).
\end{equation}
If $(-1)^n=W(\psi)$, the functional equation of $L(E/\Q,s)$ forces
$L(\chi\psi,1)=0$.  On the other hand, the main result of \cite{Ro}
shows that if $(-1)^n=-W(\psi)$, then $L(\chi\psi,1)\not=0$ for all
but finitely many $\chi$. The claim now follows from the reciprocity
law of Coates-Wiles, which relates the localization of the elliptic
units to the special value of twists of $L(\psi,s)$ (see \cite[Theorem
5.1]{RP} for example).
\end{proof}

\begin{Rem} \label{R:trivial}
The equality \eqref{E:rootno} (which is visibly incorrect when $n=0$)
is the reason for placing the trivial character in $\Xi_n^-$.  In
particular, our definitions differ from those of \cite[Definition
3.1]{RP}. 

The reader should also note that in \cite{Ru87}, the characters of
$\Gamma_n$ are indexed according to the parity of their conductors,
while we have indexed them according to the parity of their
orders. Hence our $\Xi^+$ (respectively $\Xi^-$) is denoted by $\Xi^-$
(respectively $\Xi^+$) in \cite{Ru87}.
\qed
\end{Rem}

Set $\Lambda=\co_\p[[\Gal(D_\infty/K)]]$. For any finitely generated
$\Lambda$-module $M$, we write $\cha_{\Lambda}(M)$ for the
characteristic ideal of $M$ in $\Lambda$, and $\rk_{\Lambda}(M)$ for
the $\Lambda$-rank of $M$.
Define a $\Lambda$-module
$$
X_*=\Hom_{\Z_p}(\Sel_*^\infty,\Q_p/\Z_p).
$$
Let $\iota: \Lambda \to \Lambda$ denote the canonical involution on
$\Lambda$ which is induced by inversion on group-like elements. We
adopt the convention that $\Lambda$ acts on $X_*$ via the rule
$(\lambda \cdot f)(x) = f(\lambda^{\iota} x)$ (cf. \cite[Remark
1.18]{AH}).

The following two propositions are consequences of the work of Rubin.

\begin{Prop} \label{P:mainconj}
(i) The $\Lambda$-modules $\CS_\rel$ and $X_\str$ are torsion-free of
rank one, and torsion, respectively. The $\Lambda$-module $X_{\rel}$
has rank one.

(ii) There is an equality of characteristic ideals
\begin{equation*}
\cha_{\Lambda}(X_\str) = \cha_{\Lambda} (\cS_{\rel}/\cC).
\end{equation*} 
\end{Prop}

\begin{proof} 
The fact that $\CS_\rel$ is torsion-free of the same rank as
$X_{\rel}$ may be proved exactly as in \cite[Lemma 1.1.9]{AH} (the
proof of which is essentially the same as that of \cite[Proposition
4.2.3]{PR}). The remaining claims of (i) follow from the nontriviality
of $\CC$ using the theory of Euler systems as in \cite{Ru}.  Using
(i), (ii) may be deduced from Rubin's two variable main conjecture
\cite[Theorem 4.1(ii)]{Ru91} exactly as in \cite[Proposition
2.4.16]{AH}.
\end{proof}

\begin{Prop}\label{P:local ranks}
The $\Lambda$-module $\HH^1$ is torsion free of rank $2$.  The modules
$\HH^1_\pm$ have $\Lambda$-rank $1$ and satisfy $\HH^1_+\cap\HH^1_-=0$.
The modules $H^1_\pm(D_\infty,W)$ have $\Lambda$-corank one.
\end{Prop}
\begin{proof}
Write $H^1_f(D_{\infty,\fp},W) = \varinjlim H^1_f(D_{n,\fp},W)$, where
the injective limit is taken with respect to restriction maps. Let
$\bE[\fp^{\infty}]$ and $E[\fp^{\infty}]$ denote the $\fp$-primary
torsion subgroups of $\bE$ and $E$ respectively, and write
$K(E[\fp^\infty])$ for the field obtained by adjoining the elements of
$E[\fp^{\infty}]$ to $K$. Set
$$
V_{\infty} = \Hom_{\cO_{\fp}} \left( H^1_f(D_{\infty,\fp},W),
\bE[\fp^\infty] \right),\quad
V^\pm = \Hom_{\cO_\fp} \left(
\frac{H^1_f(D_{\infty,\fp},W)}{H^1_\pm(D_{\infty,\fp}, W)},
\bE[\fp^\infty] \right).
$$
We may view $V_\infty$ and $V^\pm$ as being $\Lambda$-modules by
identifying $\Gal(D_{\infty}/K)$ with a subgroup of
$\Gal(K(E[\fp^{\infty}])/K)$ in the obvious way. It is shown in
\cite[Propositions 1.1 and 8.1]{Ru87} that the $\Lambda$-module
$V_\infty$ is torsion-free of rank $2$, while the $\Lambda$-modules
$V^\pm$ are of rank one and satisfy $V_+ \cap V_- = 0$. The
proposition now follows from the fact that fixing an identification of
$K_{\fp}/\cO_{\fp}$ with $\bE[\fp^\infty]$ induces $\Lambda$-module
isomorphisms
\begin{equation} \label{E:vtwists}
V \simeq \cH \otimes \Hom_{\cO_\fp}(\cO_\fp,T),\quad
V_\pm \simeq \cH_\pm \otimes \Hom_{\cO_\fp}(\cO_\fp,T).
\end{equation}
\end{proof}

\begin{Conj}\emph{(Rubin \cite[Conjecture 2.2]{Ru87})} \label{C:torsion}
$\cH^1 = \cH^1_+ \oplus \cH^1_-$.
\end{Conj}

\begin{Thm}\label{ranks}
We have $\rk_{\Lambda}(\cS_\pm) = \rk_{\Lambda}(X_\pm)$.
If the sign of $W(\psi)$ is $\ep$, then $X_{\ep}$
has $\Lambda$-rank one, and $X_{-\ep}$ is $\Lambda$-torsion. In
particular $\cS_{\str} = \cS_{-\ep} =0$, as $\cS_{\rel}$ is torsion-free.
\end{Thm}

\begin{proof}
Global duality (see \cite[Theorem 1.7.3]{Ru})  gives the exact sequence
\begin{equation}\label{first duality}
0\map{}\CS_\pm\map{}\CS_\rel
\map{} \HH^1/\HH^1_\pm
\map{}X_\pm\map{}X_\str\map{}0.
\end{equation}
The first claim now follows from Propositions \ref{P:mainconj} and 
\ref{P:local ranks}.

Since $\cS_{\pm} \subset \cS_\rel$, the $\Lambda$-rank of $\cS_\pm$ is
at most one, and as $\cC \subset \cS_{\ep}$ is non-trivial, we see
that the $\Lambda$-rank of $\cS_{\ep}$ is in fact equal to one. Next,
we observe that if $\rk_{\Lambda}(\cS_{-\ep}) = 1$, then
$\rk_{\Lambda}(\cS_+ \cap \cS_-) = 1$, since both $\cS_+$ and $\cS_-$
are submodules of the rank one $\Lambda$-module $\cS_\rel$. By
Proposition \ref{P:local ranks}, $\cS_+ \cap \cS_-\subset \cS_\str$,
and so also $\rk_\Lambda(\cS_\str)=1$.  But then $\cS_\rel/\cS_\str$
is a $\Lambda$-torsion module. This quotient injects into $\HH^1$
which is torsion-free by Proposition \ref{P:local ranks}. We conclude
that $\cS_\str=\cS_\rel$ and that the localization map
$\cS_\rel\map{}\HH^1$ is trivial, contradicting Proposition
\ref{P:nontriv}.

It now follows that $\rk_{\Lambda}(\cS_{-\ep}) =
\rk_{\Lambda}(\cS_{\str}) = 0$, and since $\cS_{\rel}$ is
torsion-free, this implies that both $\cS_{-\ep}$ and $\cS_{\str}$ are
equal to zero.
\end{proof}


\section{Characteristic ideals} \label{S:ideals}



\begin{Thm} \label{T:relstr}
We have the equality of characteristic ideals
$$
\cha_{\Lambda} (X_{\rel,\Lambda-\tors}) = \cha_{\Lambda}(X_\str).
$$
\end{Thm}

\begin{proof} Let $K_\infty/K$ denote the unique $\bZ_p^2$-extension
of $K$, and set $\Lambda(K_\infty):=
\cO_{\fp}[[\Gal(K_\infty/K)]]$. Write $X(K_\infty):=
\Hom(\Sel(K_\infty), \bQ_p/\bZ_p)$. It follows from Rubin's proof of
the main conjecture that $\rk_{\Lambda(K_\infty)}(X(K_\infty)) = 1$
(see e.g.  \cite[Remark 2.2]{RP} and\cite[Theorem 5.3(iii)]{Ru91}). As
$X_\rel = X$ (see \cite[Remark 3.3]{Bi}, for instance), Proposition
\ref{P:mainconj}(i) implies that $\rk_{\Lambda}(X) = 1$ also. Hence,
if $\gamma_1$ is any topological generator of
$\Gal(K_\infty/D_\infty)$, then, since
$$
\frac{X(K_\infty)}{(\gamma_1 -1)X(K_\infty)} \simeq X
$$
(see \cite[Proposition 1.2 and Theorem 2.1]{Ru85} or
\cite[p. 364--365]{Bi}), we deduce that $\gamma_1-1$ is coprime to
$\cha_{\Lambda(K_\infty)}(X(K_{\infty})_{\tors})$. The theorem now
follows directly from \cite[Theorem 3.24]{Bi} and \cite[Lemma
2.1.2]{AH} (see also \cite[Corollary 6.5]{Wi} for a more general
result along these lines). 
\end{proof}

\begin{Thm} \label{T:ideals}
Suppose that the sign of $W(\psi)$ is equal to $\ep$. Then
\begin{equation} \label{E:ideal(i)}
\cha_{\Lambda}\left(X_{\ep,\Lambda-\tors}\right)
\cha_{\Lambda}\left(\frac{\cH^{1}_{\ep}}{\cS_{\ep}}\right) =
\cha_{\Lambda}\left(\frac{\cS_\rel}{\cC}\right).
\end{equation}
If we assume that Conjecture \ref{C:torsion} is true then $\cS_\ep=\cS_\rel$
and 
\begin{equation} \label{E:ideal(ii)}
\cha_{\Lambda}\left(X_{-\ep}\right) =
\cha_{\Lambda}\left(\frac{\cH^1_\ep}{\cC} \right).
\end{equation}
\end{Thm}

\begin{proof}
From Proposition \ref{P:nontriv}, we see that $\cC \subset \cS_{\ep}$,
and so $ \rk_{\Lambda}(\cC) = \rk_{\Lambda}(\cS_\ep) = 1$.  Via global
duality (see \cite[Theorem 1.7.3]{Ru}), together with the fact that
$\cS_{\str}=0$, we have the exact sequence
\begin{equation}
0 \to \cS_{\ep} \to \cH^{1}_{\ep} \to X_\rel \to X_{\ep} \to 0.
\end{equation}
As $\cH^1/\cS_{\ep}$ is $\Lambda$-torsion, it is not hard to check
that this in turn yields the exact sequence
\begin{equation} \label{E:dagger}
0 \to \cH^{1}_{\ep}/\cS_{\ep} \to X_{\rel,\Lambda-\tors} 
\to X_{\ep,\Lambda-\tors} \to
0.
\end{equation}
The equality \eqref{E:ideal(i)} now follows
from \eqref{E:dagger} together with Theorem \ref{T:relstr} and
Proposition \ref{P:mainconj}(ii).

Now assume Conjecture \ref{C:torsion}. Then \eqref{first duality}
gives an injection $\cS_\rel/\cS_\ep\hookrightarrow \HH^1/\HH^1_\ep$
of a torsion module into a torsion-free module.  Hence $\cS_{\ep} =
\cS_\rel$.  In order to show \eqref{E:ideal(ii)}, we observe that, as
$\cS_{-\ep}=0$, \eqref{first duality} yields
\begin{equation*}
0 \to \cH^1_\ep/\cS_{\ep} \to X_{-\ep} \to X_\str
\to 0.
\end{equation*}
Combining this with the exactness of
\begin{equation*}
0 \to \cS_{\ep}/\cC \to \cH^1_\ep/ \cC \to \cH^1_\ep/\cS_{\ep} \to 0
\end{equation*}
and with Proposition \ref{P:mainconj} proves the 
equality \eqref{E:ideal(ii)}.
\end{proof}

Let $\ep$ be the sign of $W(\psi)$, and write $\overline{\psi}$ denote
the complex conjugate of the grossencharacter $\psi$. Fix a generator
$c_\ep$ of $\cH_\ep$.

\begin{Thm} \label{T:ohwell}
Assume that Conjecture \ref{C:torsion} holds. Then there exists a
generator $\cL_{-\ep}$ of $\cha_{\Lambda}(\cH^{1}_{\ep}/\cC)$ such that
the following statement is true:

Let $\chi$ be any character of $\Gamma$ of order $p^n$, where $n>0$
and satisfies $(-1)^{n+1} = W(\psi)$. Then
$$
\delta_\chi(v_\ep) \cdot\chi(\cL_{-\ep}) = \frac{L(\overline{\psi}
\chi,1)}{\Omega_E}.
$$
Here $\Omega_E \in \bR^+$ is the real period of a minimal model of
$E$, $v_\ep \in V_\ep$ is the image of $c_\ep$ under a fixed choice of
the isomorphism \eqref{E:vtwists}, and $\delta_\chi$ is the
Coates-Wiles homomorphism defined in \cite[\S2]{Ru87}. Furthermore,
$\delta_\chi(v_\ep)$ is always non-zero. 
\end{Thm}

\begin{proof} This is a direct consequence of \cite[\S10]{Ru87}, once
we fix a choice of isomorphism \eqref{E:vtwists} above. (One must also
bear in mind the last part of Remark \ref{R:trivial}.)
\end{proof}

Now Theorems \ref{T:ideals} and \ref{T:ohwell} imply that if
Conjecture \ref{C:torsion} holds, then
$$
\cL_{-\ep} \Lambda = \cha_{\Lambda}(\cH^{1}_{\ep}/\cC) =
\cha_{\Lambda}(X_{-\ep}).
$$
Hence we see that Conjecture \ref{C:torsion} implies that
$\cha_{\Lambda}(X_{-\ep})$ is generated by an element which
$p$-adically interpolates suitably normalised special values of twists
of $L(\overline{\psi},s)$, and which may therefore be viewed as being
a $p$-adic $L$-function attached to $E$.


\section{Control theorems} \label{S:control}


Define
$$
X_{n,*}=\Hom_{\co_{\p}}(\Sel_*(D_n,W),K_\p/\co_\p).
$$
Our goal in this section is to explain how to recover the
$\co_\p$-rank of $X_n$ from the $\Lambda$-modules $X_\pm$.

Fix a topological generator $\gamma\in\Gal(D_\infty/K)$ and define
$$
\omega^+_n=\prod_{1\le k\le n, k\mathrm{\ even}}
\Phi_{p^k}(\gamma)\hspace{1cm}
\omega^-_n=(\gamma-1)\prod_{1\le k\le n, k\mathrm{\ odd}}\Phi_{p^k}(\gamma)
$$
where $\Phi_{p^k}$ is the $p^k$-th cyclotomic polynomial.  Since
$\chi(\omega_n^\pm)=0$ for every $\chi\in\Xi_n^\pm$, we have
\begin{eqnarray}
\omega_n^\mp\cdot \bE(D_{n,\p})&\subset& E_\pm(D_{n,\p})
\label{first projection}\\
\omega_n^\mp \cdot \Sel(D_{n},W)&\subset& 
\Sel_\pm(D_{n},W)\nonumber
\end{eqnarray}
and similarly $\omega_n^\pm\cdot  E_\pm(D_{n,\p})=0$.

\begin{Lem}\label{local control}
The natural map 
$$
f_n: H^1_\pm(D_{n,\p},W)\map{}H^1_\pm(D_{\infty,\p},W)[\omega_n^\pm]
$$
is injective, and the $\co_\p$-corank of the cokernel of $f_n$ is a
bounded, non-decreasing function of $n$. If Conjecture \ref{C:torsion}
holds then the cokernel of $f_n$ is finite for all $n$.
\end{Lem}

\begin{proof}
Let $L$ denote the extension of $K_\fp$ obtained by adjoining
$E[\p]$ to $K_\fp$. Then it follows from
Lubin-Tate theory that $L/K_\fp$ is a totally ramified extension of
degree $p^2-1$. Hence $L \cap D_{\infty,\fp} = K_\fp$, and we 
deduce that $H^0(D_{\infty,\fp},W)=0$.
From the inflation-restriction sequence we deduce 
that 
$$
H^1(D_{n,\p},W)\map{}H^1(D_{\infty,\p},W),
$$
and therefore also $f_n$, is injective. To prove the rest
of the lemma, we compare the $\cO_\fp$-coranks of
$H^1_\pm(D_{n,\p},W)$ and $H^1_\pm(D_{\infty,\p},W)[\omega_n^\pm]$.

From Proposition \ref{P:local ranks} and the general structure theory of
$\Lambda$-modules, we see that the $\cO_\fp$-corank of
$H^1_\pm(D_{\infty,\p},W)[\omega_n^\pm]$ is equal to
$
\rk_{\cO_{\fp}} \left(\Lambda/\omega^\pm_n\Lambda \right) + e(n),
$
where $e(n)$ is a non-decreasing, bounded function of $n$. If
Conjecture \ref{C:torsion} holds, then the $\Lambda$-module
$H^1_\pm(D_\infty,W)$ is cotorsion-free, and so $e(n) = 0$ for all
$n$. On the other hand, there is an isomorphism of 
$K_\fp[\Gal(D_{n,\fp}/K_{\fp})]$-modules
\begin{equation} \label{E:galiso}
\bE(D_{n,\fp}) \otimes_{\cO_\fp} K_{\fp} \simeq D_{n,\fp}
\simeq K_\fp[\Gal(D_{n,\fp}/K_{\fp})],
\end{equation}
in which  the first isomorphism is induced by the 
logarithm of the formal group $\bE$, and
the second follows from the normal basis theorem of Galois
theory. This implies that the $\co_\p$-corank
of $H^1_\pm(D_{n,\p},W)$ is equal to the $\co_\p$-rank of 
$\Lambda/\omega_n^\pm$. The result now follows immediately.
\end{proof}

The following result is an anticyclotomic analogue of 
Kobayashi's control theorem (see \cite[Theorem 9.3]{Ko}).

\begin{Thm}
\label{K control}
The natural map
\begin{equation} \label{E:control}
X_\pm/\omega_n^\pm X_\pm\map{}X_{n,\pm}/\omega_n^\pm X_{n,\pm}
\end{equation}
is surjective. The $\co_\p$-rank of the kernel is a bounded,
nondecreasing function of $n$. If Conjecture \ref{C:torsion} holds,
then the kernel is finite for all $n$.
\end{Thm}
\begin{proof}

Write
\begin{eqnarray*}
L_{n,\pm}&=&H^1(D_{n,\p},W)[\omega_n^\pm]\big/
H^1_\pm(D_{n,\p},W),\\
L_{\infty,\pm}&=& H^1(D_{\infty,\p},W)[\omega_n^\pm]
\big/H^1_\pm(D_{\infty,\p},W)[\omega_n^\pm].
\end{eqnarray*}
Let $\mathbf{K}$ denote the maximal extension of $K$ unramified outside
$\mathfrak{f}\p$ and set
$$
H^1(D_{n,\ff},W) = \oplus_{v \mid \ff} H^1(D_{n,v},W),\quad
H^1(D_{\infty,\ff},W) = \oplus_{v \mid \ff} H^1(D_{\infty,v},W).
$$
Consider the following commutative diagram with exact rows:
$$
\xymatrix{
0 \ar[r] & {\Sel_\pm(D_n,W)[\omega_n^\pm]}\ar[r] \ar[d] & 
H^1(\mathbf{K}/D_n,W)[\omega_n^\pm] \ar[r]\ar[d] &
H^1(D_{n,\mathfrak{f}}, W)\oplus L_{n,\pm}\ar[d]\\
0 \ar[r] & {\Sel_\pm(D_\infty,W)[\omega_n^\pm]}\ar[r] & 
H^1(\mathbf{K}/D_\infty,W)[\omega_n^\pm] \ar[r] &
H^1(D_{\infty,\mathfrak{f}}, W)\oplus
L_{\infty,\pm}}
$$

The left-hand vertical arrow of this diagram is the dual of the map
\eqref{E:control}. Since $H^1(\mathbf{K}/D_n,W)\iso
H^1(\mathbf{K}/D_\infty,W)[\gamma^{p^n}-1]$, the middle vertical arrow
is an isomorphism. To prove the theorem, it therefore suffices (by the
Snake Lemma) to show that the $\cO_{\fp}$-corank of the kernel of the
right-hand arrow is a bounded, non-decreasing function of $n$, and is
finite for all $n$ if Conjecture \ref{C:torsion} holds.

For any place $v$ of $D_\infty$ dividing $\mathfrak{f}$, the extension
$D_{\infty,v}/D_{n,v}$ is either trivial (in which case there is
nothing to check) or is the unique unramified $\Z_p$-extension of
$D_{n,v}$. Assume we are in the latter case. The kernel of
$$
H^1(D_{n,v},W)\map{}H^1(D_{\infty,v},W)
$$ 
is isomorphic to $H^1(D_{\infty,v}/D_{n,v},
E(D_{\infty,v})[p^\infty])$, which is  isomorphic to a quotient
of $E(D_{\infty,v})[p^\infty]$. Since the Galois module $W$ is
ramified at all primes dividing $\mathfrak{f}$, it follows that
$E(D_{\infty,v})[p^\infty]$ is a proper $\cO_\fp$ submodule of
$W$. This implies that $E(D_{\infty,v})[p^\infty]$ is finite, because
$W$ is cofree of corank one over $\cO_\fp$, and so any proper
submodule of $W$ is finite.

In order to control the kernel of $L_{n,\pm}\map{} L_{\infty,\pm}$, we
apply the Snake Lemma to the diagram
$$
\xymatrix{
0\ar[r]&{H^1_\pm(D_{n,\p},W)}\ar[r] \ar[d] & 
H^1(D_{n,\p},W)[\omega_n^\pm] \ar[r]\ar[d] & L_{n,\pm}\ar[d]\ar[r]&0\\
0\ar[r]&{H^1_\pm(D_{\infty,\p},W)[\omega_n^\pm]}\ar[r] & 
H^1(D_{\infty,\p},W)[\omega_n^\pm] \ar[r] &L_{\infty,\pm}\ar[r]&0.}
$$
Just as in the proof of Lemma \ref{local control}, we deduce from the
inflation-restriction sequence that the middle vertical arrow of this
diagram is injective; the same inflation-restriction sequnce also
shows that that this arrow is surjective. We therefore deduce from
Lemma \ref{local control} that the $\cO_\fp$-corank of the kernel of
the right-hand vertical arrow of this diagram is a bounded,
non-decreasing function of $n$, and is finite for all $n$ if
Conjecture \ref{C:torsion} holds. This completes the proof.
\end{proof}

\begin{Prop}\label{pm control}
For any $n$, the natural map
$$
X_n\map{} (X_{n,+}/\omega_n^+ X_{n,+})\oplus (X_{n,-}/\omega_n^- X_{n,-})
$$
has finite kernel and cokernel.  
\end{Prop}
\begin{proof}
Consider the dual map
\begin{equation}\label{rec dual}
\Sel_+(D_n,W)[\omega_n^+]\oplus \Sel_-(D_n,W)[\omega_n^-]\map{} \Sel(D_n,W).
\end{equation}
By (\ref{first projection}) and the equality
$\omega_n^\pm\omega_n^\mp=\gamma^{p^k}-1$,
there is an inclusion
$$
\omega_n^-\cdot \Sel(D_n,W)+\omega_n^+\cdot \Sel(D_n,W)
\subset \Sel_+(D_n,W)[\omega_n^+]+\Sel_-(D_n,W)[\omega_n^-].
$$
Since
$$
\Sel(D_n,W)/\big(
\omega_n^-\cdot \Sel(D_n,W)+\omega_n^+\cdot \Sel(D_n,W)
\big)
$$
is a module of cofinite type over the finite ring
$\Lambda/(\omega_n^+, \omega_n^-)$, it is finite, and therefore the
same is true of the cokernel of (\ref{rec dual}).  The kernel of
\eqref{rec dual} is isomorphic to
$$
\Sel_+(D_n,W)[\omega_n^+]\cap\Sel_-(D_n,W)[\omega_n^-]
$$ 
which is again a cofinite type module over
$\Lambda/(\omega_n^+,\omega_n^-)$, and is therefore also finite.
\end{proof}

Combining Theorem \ref{ranks}, 
Theorem \ref{K control}, and Proposition \ref{pm control}
we obtain the following result.

\begin{Thm}\label{final}
Let $\ep$ be the sign of $W(\psi)$.  There is an integer $e$,
independent of $n$, such that
$$
\mathrm{corank}_{\co_\p}(\Sel(D_n,W)) = 
\mathrm{rank}_{\co_\p}(\Lambda/\omega_n^\ep)+e
$$
for $n\gg 0$. If Conjecture \ref{C:torsion} holds
then the $\cO_\fp$-corank of $\Sel(D_n,W)$ is equal to
$$
\mathrm{rank}_{\co_\p}(\Lambda/\omega_n^\ep)+
\mathrm{rank}_{\co_\p}(Y_+/\omega_n^+Y_+)+
\mathrm{rank}_{\co_\p}(Y_-/\omega_n^-Y_-)
$$
for all $n$, where $Y_\pm$ is the $\Lambda$-torsion submodule
of $X_\pm$. \qed
\end{Thm}

Theorem A of the Introduction now follows from the first part of
Theorem \ref{final}.


\begin{thebibliography}{1}

\bibitem{AH}
A. Agboola, B. Howard.
\newblock Anticyclotomic Iwasawa theory for CM elliptic curves.
\newblock Preprint (2003).

\bibitem{Bi}
P. Billot.
\newblock Quelques aspects de la descente sur une courbe elliptique
dans le cas de r\'eduction supersinguli\`ere.
\newblock {\em Comp. Math.} 58 (1986) 341--369.

\bibitem{Gr}
R.~Greenberg.
\newblock On the {B}irch and {S}winnerton-{D}yer conjecture.
\newblock {\em Invent. Math.} 72 (1983) 241--265.

\bibitem{Ka}
K. Kato.
\newblock $p$-adic Hodge theory and values of zeta functions of modular
forms.
\newblock In: Cohomologies $p$-adiques et applications arithm\'etiques III.
{\em Asterisque} 295 (2004), ix, 117--290.

\bibitem{Ko}
S. Kobayashi.
\newblock Iwasawa theory for elliptic curves at supersingular primes.
\newblock {\em Invent. Math.} 152 (2003) 1--36.

\bibitem{PR4} B. Perrin-Riou.
\newblock Arithm\'etique des courbes elliptiques \`a r\'eduction
supersinguli\`ere en $p$.
\newblock {\em Experiment. Math.} 12 (2003), no. 2, 155--186.

\bibitem{PR5} B. Perrin-Riou.
\newblock Fonctions $L$ $p$-adiques d'une courbe elliptique et points
rationnelles.
\newblock {\em Ann. Inst. Fourier} 43 (1993) 945--995.

\bibitem{PR} B. Perrin-Riou.
\newblock Th\'eorie d'Iwasawa et hauteurs $p$-adiques.
\newblock {\em Invent. Math.} 109 (1992) 137--185.

\bibitem{Po}
R. Pollack,
\newblock On the $p$-adic $L$-function of a modular form at a
supersingular prime.
\newblock {\em Duke Math. J.}  118 (2003) 523--558.


\bibitem{Ro}
D. ~Rohrlich.
\newblock On $L$-functions of elliptic curves and anticyclotomic
towers.
\newblock {\em Invent. Math.,} 75 (1984) 383--408.

\bibitem{Ru85}
K. Rubin.
\newblock Elliptic curves and $\bZ_p$-extensions.
\newblock {\em Comp. Math.} 56 (1985) 237--250.

\bibitem{Ru87}
K.~Rubin.
\newblock Local units, elliptic units, Heegner points and elliptic
curves.
\newblock {\em Invent. Math.} 88 (1987) 405--422.

\bibitem{Ru91}
K.~Rubin.
\newblock The ``main conjectures'' of
{I}wasawa theory for imaginary quadratic
  fields.
\newblock {\em Invent. Math.} 103 (1991) 25--68.

\bibitem{Ru}
K.~Rubin.
\newblock {\em Euler Systems}.
\newblock Princeton University Press, (2000).

\bibitem{RP}
K.~Rubin, R.~Pollack.
\newblock The main conjecture for elliptic curves at supersingular
primes.
\newblock {\em Ann. Math.} 159 (2004) no. 1, 447--464. 

\bibitem{Wi}
K. Wingberg.
\newblock Duality theorems for abelian varieties over
$\mathbf{Z}_p$-extensions.
\newblock {\em Advanced Studies in Pure Mathematics} 17 (1989) 471--492.

\end{thebibliography}
\end{document}